\documentstyle[11pt,figure]{article}

\begin{document}

%blackboard bold for 10 pt
%\font\bbbld=msym10
%blackboard bold for 11 pt
%\font\bbbld=msym10 scaled\magstephalf
%blackboard bold for 12 pt
\font\bbbld=msbm10 scaled\magstep1
\newcommand{\bfR}{\hbox{\bbbld R}}
\newcommand{\bfC}{\hbox{\bbbld C}}
\newcommand{\bfZ}{\hbox{\bbbld Z}}
\newcommand{\bfH}{\hbox{\bbbld H}}
\newcommand{\bfQ}{\hbox{\bbbld Q}}
\newcommand{\bfN}{\hbox{\bbbld N}}
\newcommand{\bfP}{\hbox{\bbbld P}}
\newcommand{\bfT}{\hbox{\bbbld T}}
\def\Sym{\mathop{\rm Sym}}
\newcommand{\suchthat}{\mid}
\newcommand{\halo}[1]{\Int(#1)}
\def\Int{\mathop{\rm Int}}
\def\Re{\mathop{\rm Re}}
\def\Im{\mathop{\rm Im}}
\newcommand{\union}{\cup}
\newcommand{\goesto}{\rightarrow}
\newcommand{\bdy}{\partial}
\newcommand{\n}{\noindent}
\newcommand{\p}{\hspace*{\parindent}}

\newtheorem{theorem}{Theorem}[section]
\newtheorem{assertion}{Assertion}[section]
\newtheorem{proposition}{Proposition}[section]
\newtheorem{lemma}{Lemma}[section]
\newtheorem{definition}{Definition}[section]
\newtheorem{claim}{Claim}[section]
\newtheorem{corollary}{Corollary}[section]
\newtheorem{observation}{Observation}[section]
\newtheorem{conjecture}{Conjecture}[section]
\newtheorem{question}{Question}[section]

\newbox\qedbox
\setbox\qedbox=\hbox{$\Box$}
\newenvironment{proof}{\smallskip\noindent{\bf Proof.}\hskip \labelsep}%
			{\hfill\penalty10000\copy\qedbox\par\medskip}
\newenvironment{remark}{\smallskip\noindent{\bf Remark.}\hskip \labelsep}%
			{\hfill\penalty10000\copy\qedbox\par\medskip}
\newenvironment{example}{\smallskip\noindent{\bf Example.}\hskip \labelsep}%
			{\hfill\penalty10000\copy\qedbox\par\medskip}
\newenvironment{proofspec}[1]%
		      {\smallskip\noindent{\bf Proof of Theorem 1.1.}
			\hskip \labelsep}%
			{\nobreak\hfill\hfill\nobreak\copy\qedbox\par\medskip}
\newenvironment{proofspec2}[1]%
		      {\smallskip\noindent{\bf Proof of Theorem 1.2.}
			\hskip \labelsep}%
			{\nobreak\hfill\hfill\nobreak\copy\qedbox\par\medskip}
\newenvironment{acknowledgements}{\smallskip\noindent{\bf Acknowledgements.}%
	\hskip\labelsep}{}

\setlength{\baselineskip}{1.2\baselineskip}

\title{Minimal Surfaces with Catenoid Ends}
\author{Jorgen Berglund \\Department of Mathematics
			\\University of Massachusetts
			\\Amherst, MA   01003
			\\berglund@smectos.gang.umass.edu
	 \and Wayne Rossman \\ Mathematical Institute
			\\Tohoku University
			\\Sendai 980, Japan
			\\wayne@smectos.gang.umass.edu}

\maketitle

\begin{abstract} In this paper, we use the conjugate surface construction 
to prove the existence of certain 
non-periodic symmetric immersed minimal 
surfaces.  These surfaces have finite total curvature and 
embedded catenoid ends.  Their most interesting feature is that 
they have positive genus 
yet maintain the symmetry of their genus-zero counterparts 
constructed by Jorge-Meeks and Xu.  
\end{abstract}

\section{Introduction}

In the last century, O. Bonnet, and later H. A. Schwarz, 
were the first to study 
the associate family of a minimal surface (\cite{Ni2}, \cite{Scz}).  
More recently, 
A. Schoen, H. Karcher, and others have used properties of the associate 
family to develop a method for constructing periodic
minimal surfaces (\cite{Ka1}, 
\cite{Ka2}, \cite{Ka3}, \cite{Ka4}, \cite{Kr}).  This method uses the 
particular member of the associate family known as the conjugate 
surface, and is referred to, by Karcher, as the 
{\em Conjugate Plateau Construction}.  

This construction is used to prove the results here.  
The technique begins by considering the boundary contour 
of the conjugate of a fundamental piece of the 
surface.  The contour can be described as partially unbounded 
boundary data over an unbounded convex domain.  We extend results 
of J. C. C. Nitsche \cite{Ni1} 
and Jenkins and Serrin \cite{JeSe} 
to this setting and obtain the existence of a unique minimal 
surface with this given boundary.  The existence of the original surface 
can then be argued.  

\begin{figure}
        \hspace{1.7in}
        \epsfxsize=1.6in
        \epsffile{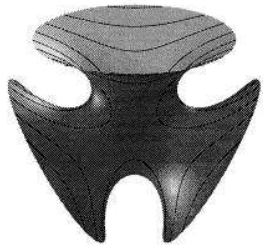}
        \hfill
%        \epsfxsize=1.6in
%        \epsffile{\psdirscr g0tritop.ps}
%        \hspace{0.9in}
	\caption{The trinoid}	
	\label{g0tri}
\end{figure}

\begin{figure}
	\hspace{1.65in}
        \epsfxsize=1.7in
        \epsffile{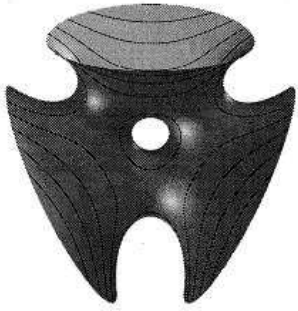}
        \hfill
%        \epsfxsize = 1.7in
%        \epsffile{\psdirscr g1tritop.ps}
%	\hspace{0.85in}
	\caption{The symmetric trinoid of genus 1}
	\label{g1tri}
\end{figure}

Our main results concern the existence of immersed finite-total-curvature 
minimal surfaces with embedded catenoid ends and genus greater than zero:
\newcounter{num}
\begin{list}%
{\arabic{num})}{\usecounter{num}\setlength{\rightmargin}{\leftmargin}}

\item For each $n \geq 3$, there exists an $n$-oid of genus 1 
that maintains all the symmetry of the 
genus-0 $n$-oid (see Figures \ref{g0tri} and \ref{g1tri}).  

\item There exist minimal surfaces based on each of the Platonic solids.  
These surfaces are of genus $f-1$ and have $v$ catenoid ends, where $f$ 
and $v$ are the number of faces and vertices of the 
corresponding Platonic solid (see Figure \ref{g7octoid}).  
\end{list}

We also prove the nonexistence of a certain symmetric 
$n$-oid of genus $n$, while indicating why 
another type might exist (see Figures \ref{cantexist} and \ref{g3tri}).  
Throughout the paper, we use Weierstrass data to draw the surfaces 
with computer graphics \cite{MESH}, and we conclude this paper by deriving 
Weierstrass data for the trinoid of genus 1.  
These data yield numerical evidence for the existence of less symmetric 
examples (see Figure \ref{fig:108}).  

The authors wish to thank Rob Kusner for many 
helpful suggestions and critical 
readings of preliminary 
drafts.  We would also like to thank: Fusheng Wei for assistance in deriving
Weierstrass data; Martin Traizet for creating Figures \ref{g3tri} 
and \ref{g7octoid}; and David Hoffman, 
Ed Thayer, and others at G.A.N.G. 
for helpful conversations and assistance with computer 
graphics.  

\section{Tools}

The following tools will be used extensively in the subsequent proofs.  

\subsection{Conjugate Surface Construction for Minimal Surfaces}

For an immersed minimal surface $M$ in $\bfR^3$ 
with finite total curvature, there exists 
a meromorphic function $g$ and a holomorphic 1-form $\eta$ defined 
on a punctured compact Riemann surface $\bar{M} \setminus 
\{p_1,p_2,...,p_\ell\}$ such that $M$ has the parametrization 
\[ \Phi(p) = Re\int_{p_0}^{p} 
\; \left( \begin{array}{c}
	(1-g^2)\eta \\
	i(1+g^2)\eta \\
	2g\eta
	\end{array}
\right) \; , \; \; \; p \in \bar{M} \setminus \{p_1,p_2,...,p_\ell\} 
\; \; \; . \] 

We refer to $\{g,\eta\}$ as the Weierstrass data for $M$, and to $\Phi$ 
as the Weierstrass representation of $M$.  
The conjugate surface $M_{conj}$ of $M$ is the minimal surface with the 
same underlying Riemann surface $\bar{M} \setminus 
\{p_1,p_2,...,p_\ell\}$, but with Weierstrass data $\{g,i\eta\}$.  
Strictly speaking, the parametrization 
$\Phi_{conj}(p)$ may only be 
well-defined on a covering of $\bar{M} \setminus 
\{p_1,p_2,...,p_\ell\}$.  

Thus we have the maps $z \rightarrow \Phi(z)$ and $z \rightarrow 
\Phi_{conj}(z)$ from simply connected domains of 
$\bar{M} \setminus \{p_1,p_2,...,p_\ell\}$ to $M$ and 
$M_{conj}$, respectively.  This induces a covering map $\phi$, the {\em 
conjugate map}, from $M_{conj}$ to $M$.  The conjugate map $\phi$ 
has the following properties:
\begin{list}%
{\arabic{num})}{\usecounter{num}\setlength{\rightmargin}{\leftmargin}}

\item $\phi$ is an isometry;

\item $\phi$ preserves the Gauss map;

\item $\phi$ maps planar principal curves in $M_{conj}$ to planar asymptotic 
curves in $M$, and maps planar asymptotic curves in $M_{conj}$ to 
planar principal curves in $M$; that is to say, $\phi$ maps 
non-straight planar geodesics to straight lines, and vice versa.  

\end{list}

It follows from the second and third properties of $\phi$ that a planar 
geodesic is mapped by $\phi$ to a line that must be perpendicular to 
the plane containing the planar geodesic.  

We note that the conjugate of 
the conjugate of $M$ is given by the Weierstrass data $\{g,-\eta\}$, 
locally giving us the original surface reflected through the origin.  

\begin{figure}
	\hspace{0.85in}
	\epsfxsize=4.0in
	\epsffile{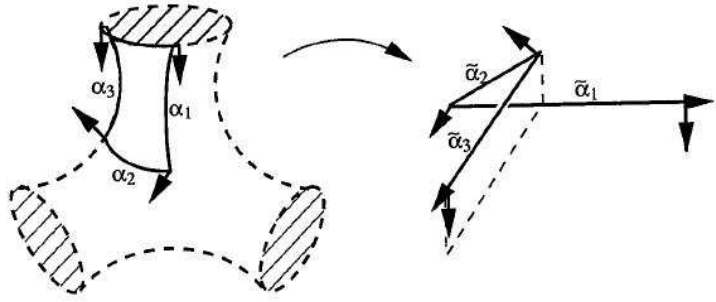}
	\hfill
	\caption{The fundamental piece of a trinoid, 
			and its conjugate surface.}
				\label{tri and conj}
\end{figure}

\begin{example}
A fundamental piece of a minimal surface is a smallest portion 
of the surface that can generate the entire surface when acted 
upon by the surface's symmetries.  
Consider a fundamental piece of an $n$-oid.  
$4n$ copies of this fundamental piece are needed to create the 
entire surface.  Note that the boundary of this fundamental piece is composed 
entirely of planar geodesics.  
The conjugate surface of this fundamental piece is thus easily 
determined by considering the properties above.  It 
is a graph over the interior of an unbounded convex 
region in some plane and thus is simply connected.  This region is 
bounded by two parallel infinite rays and one line segment connecting 
the endpoint of each ray (see Figure \ref{tri and conj}).
\end{example}

We shall say that a minimal surface has a {\em helicoid end} if 
the corresponding end of the conjugate surface is a portion of 
a catenoid end.  Thus, in the example above, the conjugate surface 
of a fundamental piece of the $n$-oid has a helicoid end.

\subsection{The Maximum Principle for Minimal Surfaces}

The following well-known lemma is the maximum principle for minimal 
surfaces.  It is a special case of a lemma by Schoen \cite{Scn}, 
and is proven there.  

\begin{lemma}
	1) {\em (Interior Maximum Principle)} Let $M_1$ and 
$M_2$ be minimal 
surfaces in $\bfR^3$.  Suppose $p$ is an interior 
point of both $M_1$, $M_2$, and suppose $T_p(M_1) = T_p(M_2)$.  
If $M_1$ lies on one side of $M_2$ near $p$, then $M_1 = M_2$.  

	2) {\em (Boundary Point Maximum Principle)} Suppose $M_1$, $M_2$ have 
$C^2$-boundaries $C_1$, $C_2$, respectively, and 
suppose $p$ is a point of both 
$C_1$, $C_2$.  Furthermore, suppose the tangent planes of both 
$M_1$, $M_2$ and $C_1$, $C_2$ agree at $p$: that is to say, suppose 
$T_p(M_1) = T_p(M_2)$, $T_p(C_1) = 
T_p(C_2)$.  If, near $p$, $M_1$ lies to one side 
of $M_2$, then $M_1 = M_2$.  
\end{lemma}

\subsection{Results on the Existence of Minimal Surfaces}

The following are two existence theorems.  
The first theorem is due to Nitsche (\cite{Ni1}, 
\cite{JeSe}), and the second 
is due to Jenkins and Serrin (\cite{JeSe}).  

\begin{theorem}
Let $D$ be a bounded convex domain in a plane.  Let $\partial \tilde{D} = 
\partial D \setminus \{p_1,...,p_r\}$.  Then there exists a solution of the 
minimal surface equation in $D$ taking on preassigned bounded continuous 
data on the arcs of $\partial \tilde{D}$.  As a surface, this solution 
contains vertical 
line segments over the jump discontinuities of the boundary data.
\end{theorem}

\begin{theorem}
{\em Monotone convergence theorem:}  
Let $\{M_n\}_{n=1}^\infty$ be a monotone 
increasing sequence of solutions of the minimal 
surface equation in a domain $D$.  If the sequence is bounded at a single 
point $p \in D$, then there exists a nonempty open set $U \subseteq D$ 
such that $\{M_n\}_{n=1}^\infty$ converges to a 
solution in $U$, and diverges to infinity on the 
complement of $U$.  The convergence is uniform on compact subsets of $U$.  
\end{theorem}

\section{Adding Handles to the $n$-oid}

The $n$-oids are well-known immersed genus-0 minimal surfaces of finite 
total curvature in $\bfR^3$ (see Figure 1 and \cite{JoMe}).  
Their defining feature is that they have $n$ 
catenoid ends, whose limiting normals span a plane $\cal P$, which
 is a plane of reflective symmetry of the surface.  
In addition, $n$-oids 
have a degree $n$ rotational symmetry about an 
axis perpendicular to $\cal P$ and a plane of reflectional symmetry 
also perpendicular to $\cal P$.  Thus the symmetry group of an 
$n$-oid is $D_n \times \bfZ_2$, the natural $\bfZ_2$-extension of 
the dihedral group.  

In this section, 
we consider the problem of adding $k$ handles to the $n$-oid, 
while preserving minimality.  We refer to these as $n$-oids of genus $k$.  

\begin{theorem}
For each $n \geq 3$, there exists an $n$-oid of genus 1 that 
maintains all the symmetries of the genus-0 $n$-oid.  
\end{theorem}

\begin{proof}
We approach the proof in the following manner: If the genus 1 surface exists, 
then it has a simply connected fundamental piece and 
the conjugate of this fundamental piece must also exist.  The boundary
contour of this conjutgate piece is among a 1-parameter family
$\hat{C}_\lambda$ of contours, each of which, we show, bounds an embedded
simply connected minimal surface.  The original fundamental piece
is, up to congruence, the conjugate to one of these, for a particular
choice of $\lambda$ which ``kills the period''.  That is to say,
$\lambda$ is chosen so that the original fundamental piece extends by 
reflection and rotation to the conjectured surface.  
We show that $n \geq 3$ is precisely 
the necessary and sufficient condition for solving this period problem.  

For any $\lambda\geq 0$ define a contour $\hat{C}_\lambda$ in $\bfR^3$ 
consisting of two straight rays and 
two line segments.  Let $\alpha_1$ be the ray $\{(0,s,0): s \geq 0\}$; 
let $\alpha_2$ be the 
line segment with endpoints at $(0,0,0)$ and $(\lambda,0,0)$; 
let $\alpha_3$ be the line segment with endpoints at $(\lambda,0,0)$ 
and $(\lambda,\cos(\frac{\pi}{n}),\sin(\frac{\pi}{n}))$; and let $\alpha_4$ 
be the ray 
$\{(\lambda+s,\cos(\frac{\pi}{n}),\sin(\frac{\pi}{n})): s \geq 0\}$;  
then $\hat{C}_\lambda$ is 
$\alpha_1 \cup \alpha_2 \cup \alpha_3 \cup \alpha_4$.  We note that the 
projection of $\hat{C}_\lambda$ to the plane $\{x_1 = 0\}$ 
lies in the boundary of an unbounded 
convex domain, which we will call $D$ (see Figure \ref{contour}).

\begin{figure}
	\hspace{1.8in}
        \epsfxsize=2.4in
%        \epsffile{\psdirhome figure_4.eps}
        \epsffile{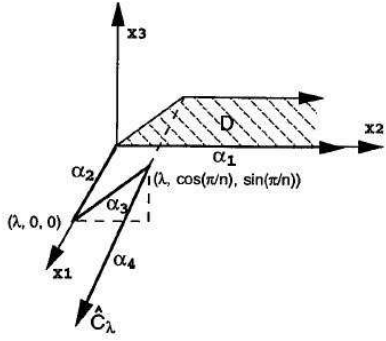}
        \hfill
	\caption{The contour $\hat{C}_\lambda$ and the convex 
		planar domain $D$.}
	\label{contour}
\end{figure}

The contour $\hat{C}_0$ is the boundary of the conjugate surface of 
the fundamental piece of the $n$-oid.  
Translating $\hat{C}_0$ in the $x_1$-direction by $(\lambda,0,0)$ we 
have a contour that we shall call $\hat{C}_{trans,\lambda}$.  (Note that 
$\hat{C}_0 = \hat{C}_{trans,0}$.)  
It follows from the known existence of the $n$-oid of genus 0 that, for all 
$\lambda$, the contour $\hat{C}_{trans,\lambda}$ bounds a minimal surface 
$\hat{M}_{trans,\lambda}$.  
The interior of $\hat{M}_{trans,\lambda}$ is a graph over the interior 
of the domain $D$.  

We shall use the surface $\hat{M}_0 = \hat{M}_{trans,0}$ to 
construct compact contours which 
converge to $\hat{C}_\lambda$.  Note that $\hat{C}_0$ coincides 
with $\hat{C}_\lambda$ along 
$\alpha_1$ and $\alpha_4$.  Choose a strictly increasing sequence 
$\lambda_i \in \bfR$ such that $\lambda_0 > \lambda$ and 
$\lim_{i \rightarrow \infty} \lambda_i = \infty$.  For each $\lambda_i$, 
choose a curve $\gamma_i$ lying in $\hat{M}_0$ with the following 
properties:
\begin{list}%
{\arabic{num})}{\usecounter{num}\setlength{\rightmargin}{\leftmargin}}

\item $\gamma_i$ has endpoints 
$(\lambda_i,\cos(\frac{\pi}{n}),\sin(\frac{\pi}{n}))$ 
and $(0,\lambda_i,0)$.  

\item $\gamma_i$ projects onto a curve $proj(\gamma_i)$ in $D$ which 
is convex with respect 
to the bounded component of $D \setminus proj(\gamma_i)$.  

\item $\gamma_i \cap \gamma_j = \phi$ for all $i \neq j$.  

\end{list}

Let $C_{\lambda_i}$ be the compact contour constructed by truncating the two 
rays $\alpha_1$, $\alpha_4$ of 
$\hat{C}_\lambda$ at $(\lambda_i,\cos(\frac{\pi}{n}),\sin(\frac{\pi}{n}))$ 
and $(0,\lambda_i,0)$, and then joining these endpoints by the 
curve $\gamma_i$.  The contour $C_{\lambda_i}$ projects to a convex 
plane curve.  
Let $D_i$ be the bounded convex region with boundary consisting of this 
projection.  $D_i \subseteq D$, and $\lim_{i \rightarrow \infty} D_i = D$ 
(see Figure \ref{contour2}).

\begin{figure}
        \hspace{1.8in}
        \epsfxsize=2.3in
        \epsffile{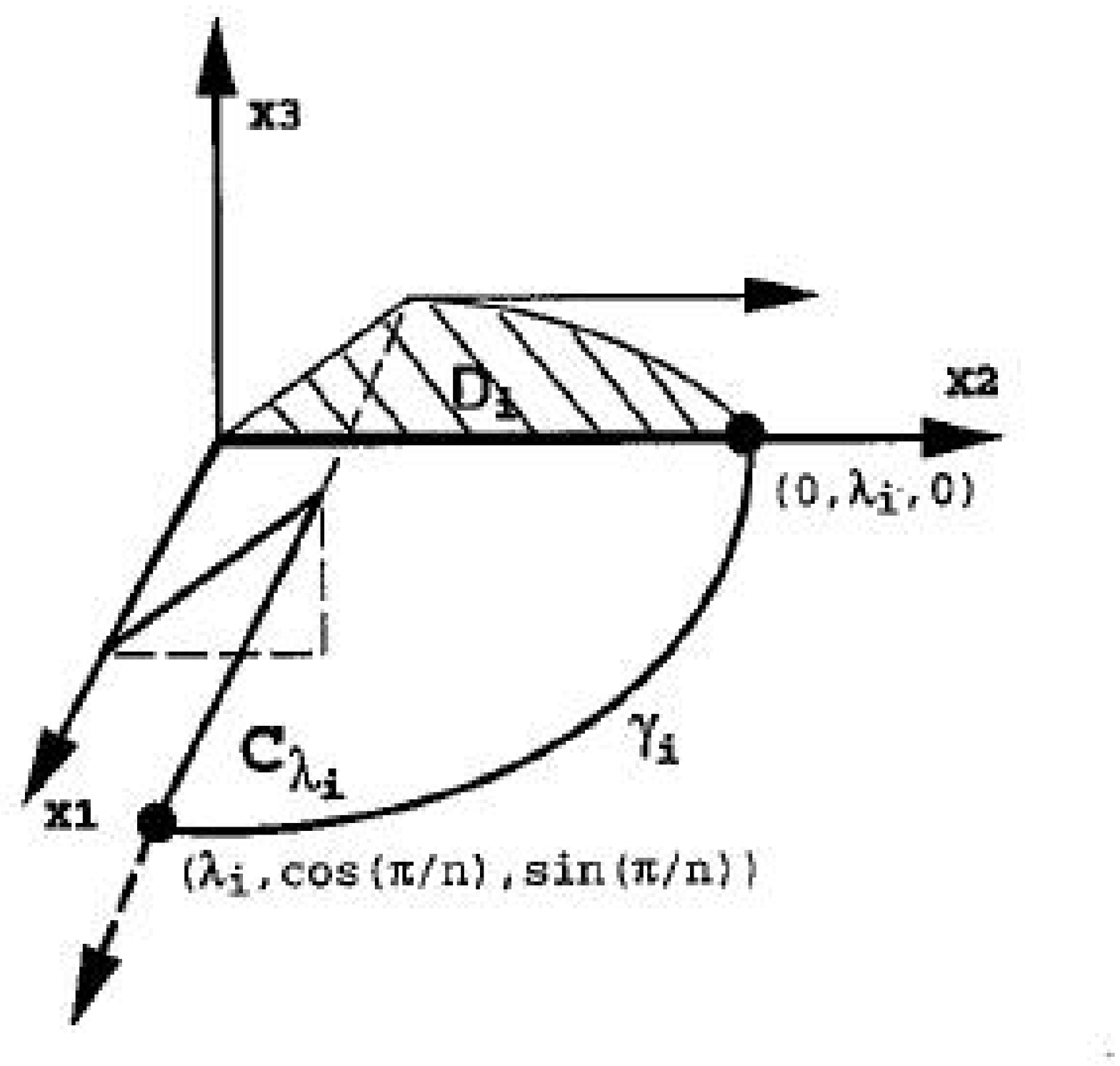}
        \hfill
 	\caption{The curve $\gamma_i$, the contour $C_{\lambda_i}$, and 
		the bounded planar region $D_i$.}	
	\label{contour2}
\end{figure}

Essentially $\{\gamma_i\}_{i=0}^\infty$ is a sequence of curves in 
$\hat{M}_0$ so that, for $j > i$, $\gamma_j$ lies ``farther out'' 
on the end of $\hat{M}_0$ than $\gamma_i$.  Since $\hat{M}_0$ 
is a graph over $D$, we conclude that $D_i \subseteq D_j$ 
for $j > i$.  Furthermore, 
since $\lambda_i \rightarrow \infty$ and $\gamma_i$ projects to a convex 
curve in $D$, 
we may conclude that $\gamma_i$ lies outside 
any given compact region in $\bfR^3$, for $i$ chosen large enough.  
Thus, $\lim_{i \rightarrow \infty} C_{\lambda_i} = \hat{C}_\lambda$.  

Since $C_{\lambda_i}$ can be viewed as piece-wise continuous boundary data 
on the domain $D_i$, discontinuous at only two points, we can apply 
Theorem 2.1.  Therefore, for each $i$, there exists a 
minimal surface bounded by $C_{\lambda_i}$.  Applying the maximum 
principle, we see that these solutions are unique, and we call these 
unique solutions $M_{\lambda_i}$.  For any fixed $i$, the surface 
$M_{\lambda_j}$, for $j \geq i$, 
restricts to a solution of the minimal surface 
equation over $D_i$.  Note that the restriction of 
$M_{\lambda_j}$ to $D_i$ 
may have different boundary data than $M_{\lambda_i}$ over $D_i$.  

{\em claim:}  Fixing a positive integer $i$, the restrictions of 
$M_{\lambda_j}$ to $D_i$, 
for $j \geq i$, form a monotonically increasing 
sequence of solutions of the minimal surface equation over $D_i$.  

{\em proof of claim:}  Let $\epsilon$ be a fixed positive number.  
By ``sliding $M_{\lambda_k}$ underneath 
$M_{\lambda_j}$'', for $i \leq k < j$, we mean this: 
We start with copies of $M_{\lambda_k}$ translated by the 
vectors $s \cdot (-1,0,+\epsilon)$ 
for $s \geq 0$, and we call these copies $M_{\lambda_k,s}$.  
Since $M_{\lambda_k}$ and $M_{\lambda_j}$ are both graphs over $D$,
it is clear that for $s > \frac{1}{\epsilon}\sin(\frac{\pi}{n})$, 
$M_{\lambda_k,s}$ and $M_{\lambda_j}$ are disjoint.  
Choose $s > \frac{1}{\epsilon}\sin(\frac{\pi}{n})$.  We then 
lower the value of $s$ until we reach the first value of $s$ so 
that $M_{\lambda_k,s} \cap M_{\lambda_j} \neq \phi$ 
(see Figure \ref{slide}).  Let 
$s_0 = \sup\{s \geq 0 \, | \, M_{\lambda_k,s} \cap M_{\lambda_j} \neq 
\phi\}$.  Proving the claim is equivalent to showing that 
$s_0 = 0$.  
Thus, we are sliding one surface underneath the other with respect to 
the positive $x_1$ direction.  

In all subsequent references to height, 
we mean height with respect to the positive $x_1$ direction.  
(In Figure \ref{slide}, we see $M_{\lambda_j}$ and a translated 
copy of $M_{\lambda_i}$.  From the point of view of the positive 
$x_3$-axis, the copy of $M_{\lambda_i}$ lies above $M_{\lambda_j}$.  
But with respect to the positive $x_1$-axis, the copy of 
$M_{\lambda_i}$ actually lies below $M_{\lambda_j}$.  We take 
the latter perspective here.)

By the interior maximum principle, this first contact 
between $M_{\lambda_k,s}$ and $M_{\lambda_j}$ occurs along 
the boundary of $M_{\lambda_k,s}$.  
The first point of contact cannot 
occur at a point of the translated copy of $\gamma_k$.  We can see this 
by showing that the original copy of $\gamma_k$ must lie 
strictly to one side of $M_{\lambda_j}$, except at its endpoints.  
The curve $\gamma_k$ is, by construction, a curve on $\hat{M}_0$.  
Again using the interior 
maximum principle, the surface $\hat{M}_0$ can be ``slid 
underneath'' $M_{\lambda_j}$ with no first point of 
contact in the interior.  In particular, 
there is no first point of contact on the interior of $\gamma_k$, and 
we conclude $\gamma_k$ lies underneath 
$M_{\lambda_j}$.  Thus the contact between 
$M_{\lambda_k,s_0}$ and $M_{\lambda_j}$ occurs in 
$\partial M_{\lambda_k,s_0}$, but cannot occur within the 
interior of the translated copy of $\gamma_k$ in 
$\partial M_{\lambda_k,s_0}$.  
Hence, the contact between $M_{\lambda_j}$ and 
$M_{\lambda_k,s_0}$ must occur along the common straight line 
boundaries of the two surfaces, and so $s_0$ must be zero.  Therefore 
$M_{\lambda_k}$ lies entirely underneath $M_{\lambda_j}$.  
This proves the claim.  

\begin{figure}
        \hspace{1.6in}
        \epsfxsize=2.7in
        \epsffile{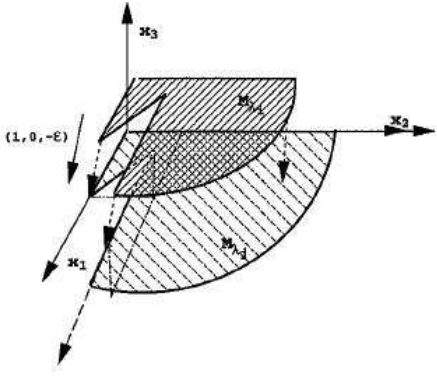}
        \hfill
 	\caption{The surface $M_{\lambda_k}$ sliding ``underneath'' 
			$M_{\lambda_j}$.}	
	\label{slide}
\end{figure}

In the proof of the claim above, if we replace 
``$M_{\lambda_k}$'' by ``$M_{\lambda_j}$'', and we replace 
``$M_{\lambda_j}$'' by ``$\hat{M}_{trans,\lambda}$'', 
we have a proof that $M_{\lambda_j}$, for any $j$, 
lies underneath $\hat{M}_{trans,\lambda}$.  Thus, over any point $p \in D_i$, 
the height of the solutions $M_{\lambda_j}$ are 
uniformly bounded above by the height of 
$\hat{M}_\lambda$ over $p$.  This allows us to apply the monotone convergence 
theorem.  

By Theorem 2.2, since the solutions $M_{\lambda_j}$, $j \geq i$, are 
monotonically increasing and are uniformly bounded-above 
over each $p \in D_i$, we conclude that the sequence of solutions 
$\{M_{\lambda_j}\}_{j=i}^\infty$, restricted to the domain 
$D_i$, converges uniformly to some solution over $D_i$.  
Since the choice of $i$ is arbitrary and 
$\lim_{i \rightarrow \infty} D_{\lambda_i} = D$, we see that the 
sequence of solutions $\{M_{\lambda_j}\}_{j=1}^\infty$ converges to 
a solution over $D$, and we call this solution $\hat{M}_\lambda$.  
Considered as a surface, $\hat{M}_\lambda$ has boundary $\hat{C}_\lambda$, 
and is a graph over $D$.  

Let $\tilde{M}_\lambda$ be the conjugate surface of $\hat{M}_\lambda$.  
Denote the 
boundary of $\tilde{M}_\lambda$ by $\tilde{C}_\lambda$, and denote each 
planar geodesic in $\tilde{C}_\lambda$ by $\tilde{\alpha}_i$, in 
correspondence with its preimage line segment or ray 
$\alpha_i \subseteq \hat{C}_\lambda$  (See Figure \ref{conjg1}).

Note that the surfaces $\hat{M}_{trans,\lambda}$ and 
$\hat{M}_0$ both have helicoid 
ends, which are asymptotic to each other.  Since $\hat{M}_\lambda$ 
lies between these two surfaces, the 
asymptotic behavior of the end of $\hat{M}_\lambda$ is determined.  
We now argue that $\hat{M}_\lambda$ has a helicoid end.  
We do this by showing that $\hat{M}_\lambda$ has finite 
total curvature (thus $\tilde{M}_\lambda$ has finite
total curvature), and then applying Schoen's result on 
complete finite-total-curvature ends \cite{Scn}.  

Choose an orientation on $M_{\lambda_i}$, and consider the Gauss map 
$G: M_{\lambda_i} \rightarrow S^2$.  Since $M_{\lambda_i}$ is a 
graph, the image Im($M_{\lambda_i}$) $\subseteq S^2$ 
of $M_{\lambda_i}$ under 
$G$ must be contained in a hemisphere.  The image 
Im($C_{\lambda_i}$) of $C_{\lambda_i} = \partial M_{\lambda_i}$ 
under $G$ is a set of curves in $S^2$.  $M_{\lambda_i}$ is a
compact surface, hence has finite total curvature.  Therefore 
$G$ is a branched covering map from $M_{\lambda_i}$ to 
Im($M_{\lambda_i}$) with finite degree.  

Let $P$ be the plane containing the points $(0,0,0)$, 
$(\lambda,0,0)$, and 
$(\lambda,\cos(\frac{\pi}{n}),\sin(\frac{\pi}{n}))$.  $M_{\lambda_i}$ 
lies to one side of $P$ at 
$(\lambda,\cos(\frac{\pi}{n}),\sin(\frac{\pi}{n}))$, thus $G$ 
cannot be branched at this point.  Furthermore, from the geometry 
of $C_{\lambda_i}$, we see that the preimage set of 
$G( \, (\lambda,\cos(\frac{\pi}{n}),\sin(\frac{\pi}{n})) \, )$ 
in $M_{\lambda_i}$ consists 
only of the point $(\lambda,\cos(\frac{\pi}{n}),\sin(\frac{\pi}{n}))$.  
Thus the degree of the covering map $G$ must be 1.  
It follows that 
the total area of Im($M_{\lambda_i}$) must be less than $2\pi$ 
for all $i$, even when the area is counted with multiplicity.  
(In fact, the area is close to $\pi$ for large values of $\lambda_i$.)  
Therefore the total curvature of $M_{\lambda_i}$ is less than 
$2\pi$ for all $i$, and the limit surface $\hat{M}_\lambda$ has 
total curvature at most $2\pi$.  In particular, $\hat{M}_\lambda$ has
finite total curvature.  

Since conjugation is an isometry, we know that 
$\tilde{M}_\lambda$ also has finite total curvature.  First we extend 
$\tilde{M}_\lambda$ by reflection through the plane containing 
$\tilde{\alpha}_1$, and then we extend further by reflection through 
the plane 
containing $\tilde{\alpha}_4$.  The resulting surface is an annulus 
with one complete end.  This end has finite total curvature.  
Thus we can apply the result of Schoen \cite{Scn} to conclude that this 
end must be either a planar end or a catenoid end.  Clearly 
$\hat{M}_\lambda$ does not have a planar end, since 
the rays $\alpha_1$ and $\alpha_4$ do not lie in a common plane.  We 
conclude that $\hat{M}_\lambda$ has a helicoid end.  

We now show that $\hat{M}_\lambda$ is the unique minimal surface with a 
helicoid end and boundary $\hat{C}_\lambda$, that is a graph over $D$.  
Assume there is another such surface $S$.  
Consider sliding $S$ underneath $\hat{M}_\lambda$; that is, 
consider the proof of the last claim, but with 
``$M_{\lambda_k}$'' replaced by ``$S$'' and ``$M_{\lambda_j}$'' 
replaced by ``$\hat{M}_\lambda$''.  
If $S$ is slid underneath $\hat{M}_\lambda$ in this way, 
contact ``at infinity'' (meaning contact at the ends) between 
 $\hat{M}_\lambda$ and a copy of $S$ 
cannot occur before the boundaries coincide.  
This follows, since by assumption $S$ has a helicoid end, 
and thus the ends of $S$ and $\hat{M}_\lambda$ are asymptotic to each other.  
Also, by the interior maximum principle, contact at a finite point cannot 
occur before the boundaries coincide.  So first contact occurs exactly 
when the boundaries coincide, and therefore $S$ lies 
underneath $\hat{M}_\lambda$.  Similarly, we can show that 
$\hat{M}_\lambda$ lies 
underneath $S$, and we conclude $S = \hat{M}_\lambda$.  

\begin{figure}
        \hspace{1.0in}
        \epsfxsize=4.0in
        \epsffile{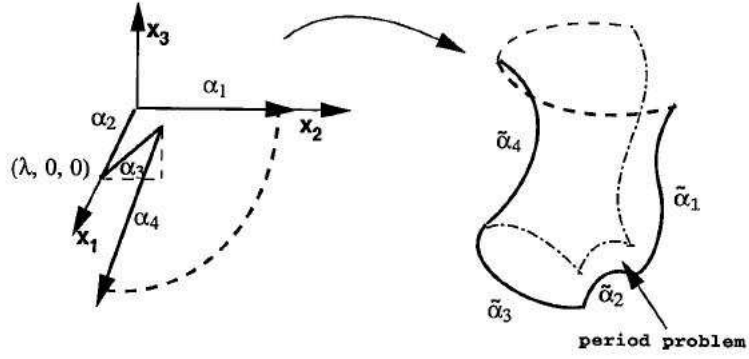}
        \hfill
 	\caption{The surface $\hat{M}_\lambda$ and its conjugate surface 
$\tilde{M}_\lambda$.}
				\label{conjg1}
\end{figure}

Since $\tilde{M}_\lambda$ is the conjugate surface of a graph over a convex 
domain, it is also a graph \cite{Kr}, and is therefore embedded.  
Also, since $\tilde{M}_\lambda$ is the conjugate of a surface with a 
single end that is a 90 degree arc of a helicoid end, 
$\tilde{M}_\lambda$ itself has a single end that is a 90 degree arc of a 
catenoid end.  If $\tilde{\alpha}_2$ and 
$\tilde{\alpha}_4$ were to lie in the same plane, then 
$\tilde{M}_\lambda$ could be extended by Schwarz reflection to a 
complete embedded minimal surface with catenoid ends.  
(The Schwarz reflection principle states that a 
minimal surface, which meets a plane orthogonally, may be extended to a 
larger minimal surface by reflection in that plane \cite{Scz}, 
\cite{HoMe}.)  This extended surface would be 
an $n$-oid with a handle symmetrically placed in the middle, 
i.e. an $n$-oid of genus 1.  However, 
$\tilde{\alpha}_2$ and $\tilde{\alpha}_4$ do not necessarily lie in the 
same plane.  This period problem can be viewed clearly by considering 
$\tilde{\alpha}_3$ lying in a plane (see Figure \ref{alpha3}).

\begin{figure}
        \hspace{1.1in}
        \epsfxsize=3.8in
        \epsffile{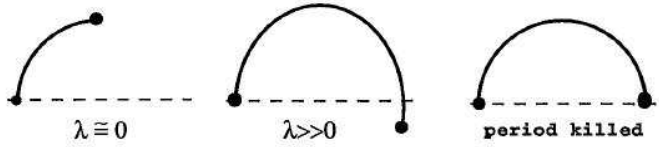}
        \hfill
 	\caption{The planar geodesic $\tilde{\alpha}_3$ lying in a plane: 
	     for $\lambda$ close to zero; for $\lambda$ large; and for 
	     $\lambda$ adjusted to kill the period}  
	\label{alpha3}
\end{figure}

As $\lambda \rightarrow 0$, the length of $\tilde{\alpha}_2$ approaches zero 
and $\tilde{\alpha}_3$ approaches a curve on the fundamental piece of the 
$n$-oid, as in the first part of Figure \ref{alpha3}.
We wish to show that for some large value 
of $\lambda$, $\tilde{\alpha}_3$ appears as in the second part of 
Figure \ref{alpha3}.
Then, 
by the Intermediate Value theorem, there will exist a value of $\lambda$ for 
which $\tilde{\alpha}_3$ appears as in the third part of 
Figure \ref{alpha3}.  
Therefore 
the period problem will be solved.

We accomplish this with a helicoidal comparison argument.  Consider a 
half-turn of a helicoid slid on ``top'' (again with respect to $x_1$ as 
height) of $\hat{M}_\lambda$, so that they share the boundary 
$\alpha_2$, $\alpha_3$, and 
$\alpha_4$, and so that $\hat{M}_\lambda$ and this half-turn of a 
helicoid lie on 
the same side of the plane $P$ through the points $(0,0,0)$, 
$(\lambda,0,0)$, 
and $(\lambda,\cos(\frac{\pi}{n}),\sin(\frac{\pi}{n}))$.  If $\lambda$ is 
sufficiently large, this ``sliding'' can be done so that the first moment 
of contact occurs along the boundary curves of the two surfaces.  
We note that the normal vectors of the helicoid and $\hat{M}_\lambda$ 
coincide at $(\lambda,0,0)$ and 
$(\lambda,\cos(\frac{\pi}{n}),\sin(\frac{\pi}{n}))$, 
but not, by the boundary point maximum principle, 
at any point of the interior of $\alpha_3$.  
Since $\hat{M}_\lambda$ lies below the helicoid, 
the normal of $M_\lambda$ must turn faster than the normal 
of the helicoid, when moving from 
$(\lambda,\cos(\frac{\pi}{n}),\sin(\frac{\pi}{n}))$ 
to $(\lambda,0,0)$ along 
$\alpha_3$.  The same is then true along the corresponding curves 
in $\tilde{M}_\lambda$ and the conjugate surface 
to the helicoid.  Recall that the catenoid is the conjugate surface to 
the helicoid, thus, the conjugate of $\alpha_3$, as a curve on the 
helicoid, is a half-circle.  Since the conjugate map is an isometry, 
the conjugate of $\alpha_3$ as a curve on $\hat{M}_\lambda$, 
and as a curve on the 
helicoid, must be of equal length.  
The curve $\tilde{\alpha}_3$ is thus forced to lie 
below the half-circle in the catenoid (see Figure \ref{compare}).  
\begin{figure}
        \hspace{1.8in}
        \epsfxsize=2.3in
        \epsffile{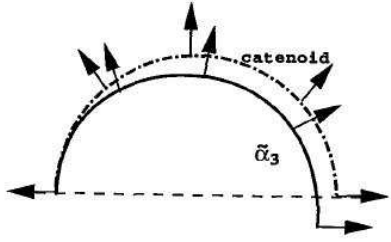}
        \hfill
	\caption{The planar geodesic $\tilde{\alpha}_3$ and a half-circle 
	     lying in the catenoid}  
	\label{compare}
\end{figure}
\end{proof}

\begin{remark}
In order to fully justify the picture we have given of the genus-1 
symmetric $n$-oid, we will show that the only branch points 
of this surface are simple branch points on curves identified, under 
reflection, with the curve $\tilde{\alpha}_1$.  

First note that the Gauss map must be $n$-to-1 on the genus-1 $n$-oid.  This 
follows from a result of Jorge-Meeks \cite{JoMe}: 
if the $n$ ends of a minimal surface are each separately embedded, then 
the degree of the Gauss map is $g + n - 1$, where $g$ is the 
genus of the underlying Riemann surface.  

It is shown in the proof above that $\lambda$ can be chosen so that 
$\tilde{M}_\lambda$ has no period problem.  Choose this value 
for $\lambda$.  
Note that $\hat{M}_\lambda$ has the same normal vector at the two extremes 
of $\alpha_1$: one of these extremes is the endpoint $(0,0,0)$; 
and the other extreme is the limit as one travels 
out to $\infty$ along the positive $x_2$-axis.  It follows that the Gauss map 
must turn back on itself at some point on $\alpha_1$, hence the Gauss 
map has a branch point at some point on $\alpha_1$.  This fact enables us 
to locate 2$n$ branch points of the Gauss map 
on the resulting symmetric $n$-oid of genus 1.  
From the Riemann-Hurwitz formula, we see that these 
are the only branch points and they must all be simple.  
\end{remark}

\begin{remark}
 Theorem 3.1 will not hold in the case $n=2$.  
In this case $\alpha_3$ is parallel to the $x_3$-axis and therefore
the half turn of the helicoid can be slid underneath $\hat{M}_\lambda$ for all
positive values of $\lambda$.  Then by helicoidal comparison we will always
have a situation comparable to the one pictured in figure \ref{compare}, except
that $\tilde{\alpha}_3$ will sit above the half circle of the catenoid,
thus the period cannot be killed.  Of course, 
Schoen \cite{Scn} has shown that the only immersed minimal surfaces 
with two catenoid ends (and no other ends) are the catenoids themselves.  
\end{remark}

A natural question to ask is whether one can add more handles to the $n$-oid, 
especially while preserving symmetry.  Our next result shows that a certain 
example is impossible.  

\begin{theorem}
There does not exist a symmetric $n$-iod of genus-$n$ such that one handle 
is situated on each of the rays originating at the center of the $n$-iod 
and directed through the center of each of the catenoid ends.  
\end{theorem}

\begin{proof}
Assume the surface exists, and let $M$ be the conjugate of 
a fundamental piece 
of the surface.  Let $\tilde{M}$ be the conjugate surface of $M$ (see 
Figure \ref{cantexist}).  Regardless of the length of $\alpha_3$ in $M$, 
we can slide a piece of a helicoid underneath $M$ (again w.r.t. 
the $x_1$ direction), so that: 
the boundary of $M$ and the boundary of the 
helicoid-piece coincide along $\alpha_1 \cup 
\alpha_2 \cup \alpha_3$; $M$ and the helicoid-piece lie on the same side 
of the plane which contains $\alpha_1 \cup 
\alpha_2 \cup \alpha_3$; and the interior of $M$ lies strictly 
to one side of the interior of the helicoid-piece.  
Then moving from the origin 
along $\alpha_2$, the normal of $M$ must rotate ahead of the 
normal of the helicoid, forcing the 
corresponding curves in the conjugate surfaces to always appear 
as in Figure \ref{compare}.  Therefore the period cannot be killed.  
\begin{figure}
        \hspace{1.03in}
        \epsfxsize=4.0in
        \epsffile{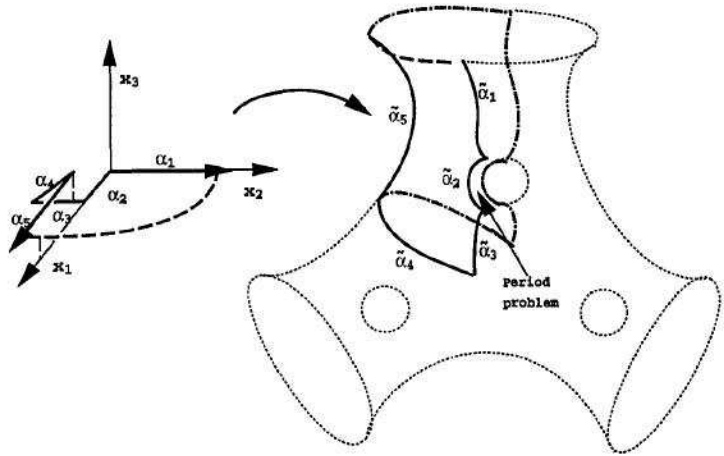}
        \hfill
	\caption{The surface $M$ and its conjugate surface $\tilde{M}$}
	\label{cantexist}
\end{figure}

\end{proof}

\begin{figure}
        \hspace{1.85in}
        \epsfxsize=2.0in
        \epsffile{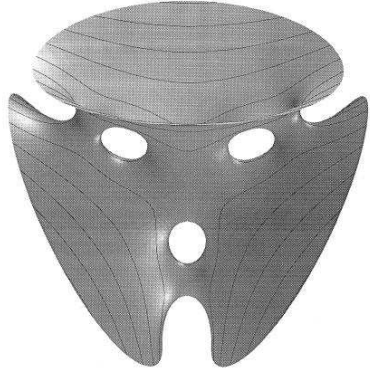}
        \hfill
        \caption{The trinoid of genus 3}		
	\label{g3tri}
\end{figure}

\begin{remark}
Using Karcher's view of the $n$-oid as the limit of a deformation
of the 2$n$-winged Scherk's towers (\cite{Ka1}, \cite{Ka3}), Martin Traizet
was able to numerically argue the existence of a symmetric $n$-oid of 
genus $n$ with handles situated on rays originating from the center and
bisecting the axes of the catenoid ends \cite{Tr} (see Figure \ref{g3tri}).
\end{remark}

\section{Minimal Surfaces Based on the Platonic Solids}

\begin{figure}
        \hspace{0.34in}
        \epsfxsize=2.0in
        \epsffile{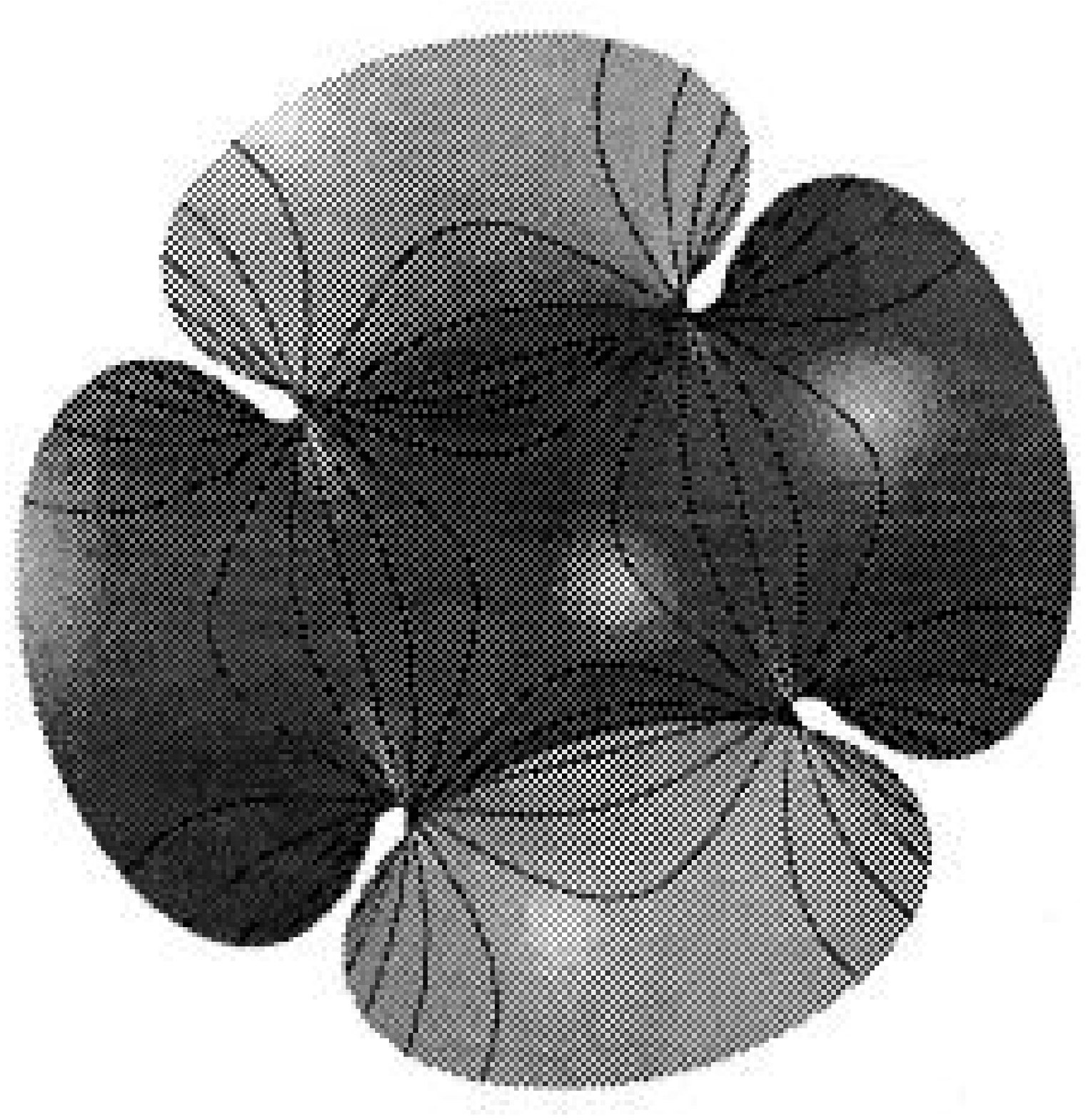}
        \hfill
        \epsfxsize=2.0in
        \epsffile{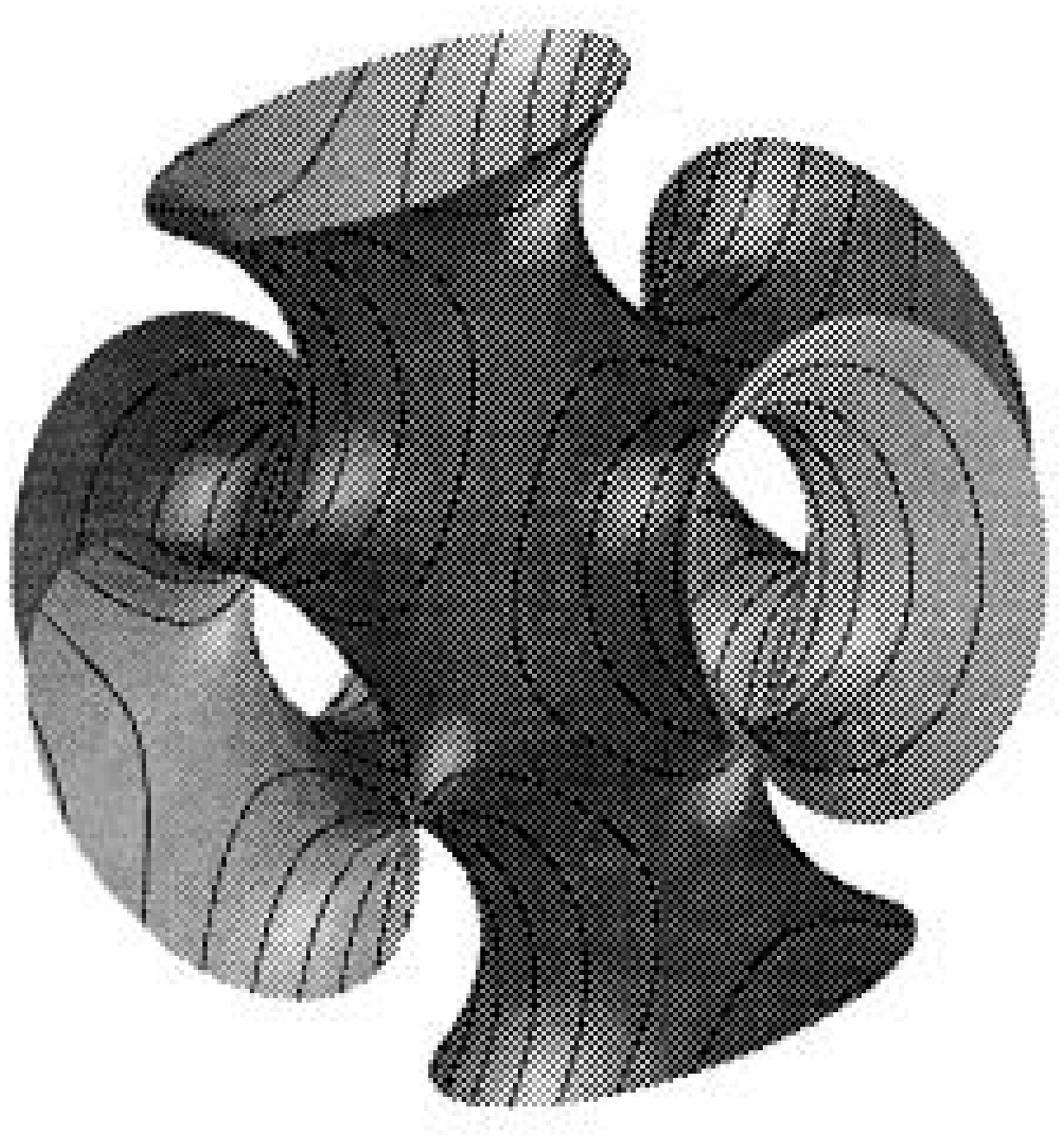}
        \hspace{0.46in}
	\caption{The tetroid and the octoid}
	\label{octoid}
\end{figure}

The method of the proof of Theorem 3.1 can be used to prove the existence 
of other immersed finite total curvature surfaces with catenoid ends 
and genus greater than zero.  Recently Y. Xu \cite{Xu} constructed
genus-0 immersed minimal surfaces based on the 
Platonic solids.  Topologically, they 
can be thought of as the surface of each Platonic solid with a catenoid end 
replacing each vertex of the solid (see Figure \ref{octoid}).
As before, we rely on the existence of the 
genus zero surfaces to construct their higher genus counterparts.

Let $\cal T$, $\cal C$, $\cal O$, $\cal D$ and 
$\cal I$ be the symmetry 
groups of the Platonic solids: the tetrahedron, cube, octahedron,
dodecahedron, and 
icosahedron, respectively.  Note that $\cal C$ and 
$\cal O$ are isomorphic, as are $\cal D$ and $\cal I$.

\begin{figure}
        \hspace{1.62in}
        \epsfxsize=2.5in
        \epsffile{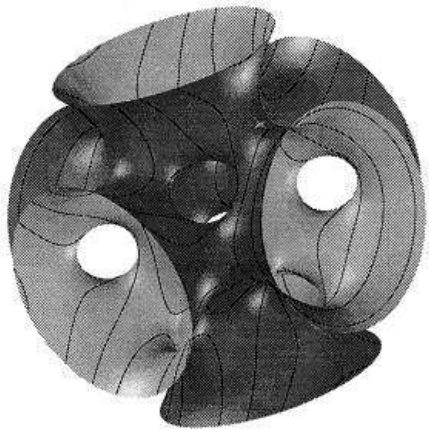}
        \hfill
 	\caption{The octoid of genus 7}		
	\label{g7octoid}
\end{figure}

\begin{theorem}
The following minimal surfaces with catenoid 
ends and finite total curvature exist: 
\begin{list}%
{\arabic{num})}{\usecounter{num}\setlength{\rightmargin}{\leftmargin}}

\item A genus-3 surface with 4 ends and symmetry group isomorphic to 
$\cal T$.  

\item A genus-5 surface with 8 ends and symmetry group isomorphic to 
$\cal C$.  

\item A genus-7 surface with 6 ends and symmetry group isomorphic to 
$\cal O$ {\em (see Figure \ref{g7octoid})}.  

\item A genus-11 surface with 20 ends and symmetry group isomorphic to 
$\cal D$.  

\item A genus-19 surface with 12 ends and symmetry group isomorphic to 
$\cal I$.  
\end{list}
\end{theorem}

\begin{proof}
Quotienting by all symmetries, we have the fundamental piece, and its 
conjugate, of the genus-0 examples.  They appear 
in Figure \ref{platconj}, where 
$\theta_1$, $\theta_2$ depend on which Platonic solid we consider.  
Consider the contour $C$ and its conjugate (see Figure \ref{gplatconj}).  When 
projected onto a plane perpendicular to $\alpha_2$ and $\alpha_4$, the 
contour projects to part of the boundary of an unbounded convex domain.  
Viewing the direction of 
$\alpha_4$ as the height, the existence and uniqueness of a minimal graph with 
boundary $C$ and a helicoid end follow 
as in the proof of Theorem 3.1.  Looking at the conjugate of this 
surface, we see that we have a single period problem.  Again, since 
$\theta_1 < \frac{\pi}{2}$, the helicoidal 
comparison argument shows that this period problem can be solved.  

\begin{figure} 
        \hspace{0.9in}
        \epsfxsize=4.0in
        \epsffile{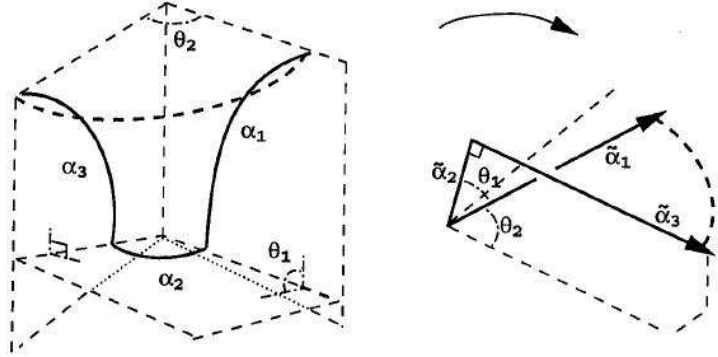}
        \hfill
	\caption{A fundamental piece of the genus-0 surface, and its conjugate 
	     surface}		
	\label{platconj}
\end{figure}

\begin{figure}
        \hspace{0.9in}
	\epsfxsize=4.0in
        \epsffile{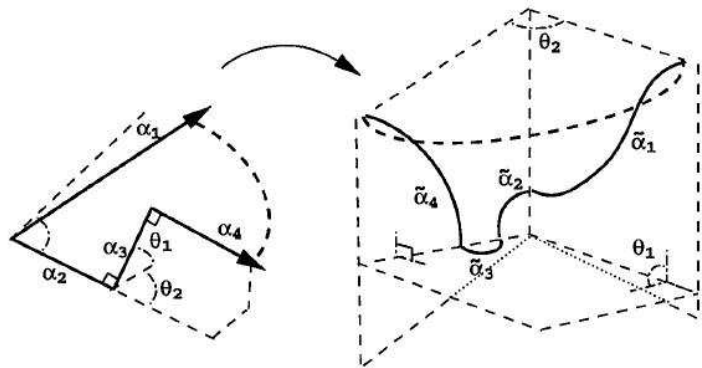}
        \hfill
	\caption{The fundamental piece with 
		contour $C = \alpha_1 \cup \alpha_2 \cup \alpha_3 
			\cup \alpha_4$, and its conjugate, which is a 
		fundamental piece of the higher genus surface}  
	\label{gplatconj}
\end{figure}

We list below the values of the angles $\theta_1$ and $\theta_2$ for each 
of the five surfaces:  
\begin{list}%
{\arabic{num})}{\usecounter{num}\setlength{\rightmargin}{\leftmargin}}

\item tetroid 	$\theta_1 = \frac{\pi}{3}$, $\theta_2 = \frac{\pi}{3}$.
\item cuboid 	$\theta_1 = \frac{\pi}{3}$, $\theta_2 = \frac{\pi}{4}$.
\item octoid 	$\theta_1 = \frac{\pi}{4}$, $\theta_2 = \frac{\pi}{3}$.
\item dodecoid 	$\theta_1 = \frac{\pi}{3}$, $\theta_2 = \frac{\pi}{5}$.
\item icosoid 	$\theta_1 = \frac{\pi}{5}$, $\theta_2 = \frac{\pi}{3}$.
\end{list}
\end{proof}

\begin{remark}
Xu also constructs minimal surfaces which can be thought of as Platonic
solids with catenoid ends added to the edges.  While in these cases
we can show that
the contour $C$ bounds a minimal surface, we find we can not apply the
helicoidal comparison test to kill the period 
as the angle between $\alpha_2$ and $\alpha_3$,
as well as the angle between $\alpha_3$ and $\alpha_4$, is not $\frac{\pi}{2}$.
\end{remark}

\section{Weierstrass Data for the Genus-1 Trinoid}

To construct Weierstrass data for the genus-1 trinoid, 
it is convenient to consider one fourth of 
the surface obtained by quotienting by the reflection in $\cal P$ and 
the reflection in one other reflectional plane of symmetry.  

\begin{figure}
        \hspace{0.92in}
        \epsfxsize=4.0in
        \epsffile{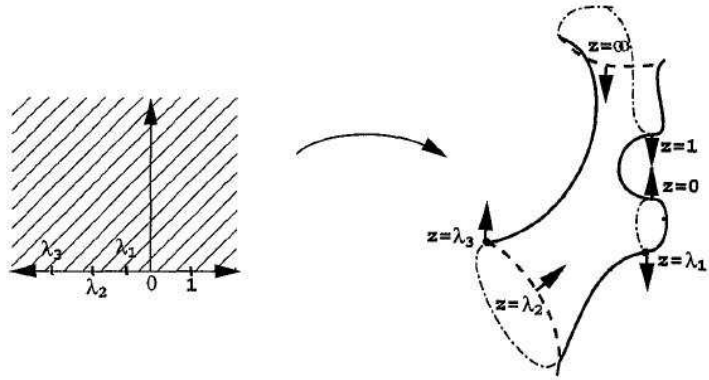}
        \hfill
	\caption{The image of the Weierstrass Data in equation (5.1)}
	\label{w-map}
\end{figure}

Let $\bar{M}$ be $\{z \in \bfC: Im(z) \geq 0\}$.  
Consider the Weierstrass data satisfying 
\begin{equation}
g^2 = c \frac{(z-1)(z-\lambda_1)}{z(z-\lambda_3)^2} \; ,  \; \; \; \; 
\eta = \frac{(z-\lambda_3)dz}{g(z-\lambda_2)^2} \; ,
\end{equation}
with $0 > \lambda_1 > \lambda_2 > \lambda_3$.  With 
this Weierstrass data, integrating over $\bar{M}$, 
we have the surface that is one fourth of the genus-1 
trinoid (not necessarily fully symmetric), 
up to some period problems.  The constant $c$ 
can be chosen so that the angle between the normals 
at the ends of the surface is $\frac{2\pi}{3}$ (see Figure \ref{w-map}).
The constant $c$ is a positive real, and its exact 
value is 
\[ c = \frac{-3\lambda_2(\lambda_2-\lambda_3)^2}
{(\lambda_2-1)(\lambda_2-\lambda_1)} \; \; . \]

\begin{figure}
	\hspace{0.5in}
	\epsfxsize = 2.0in
	\epsffile{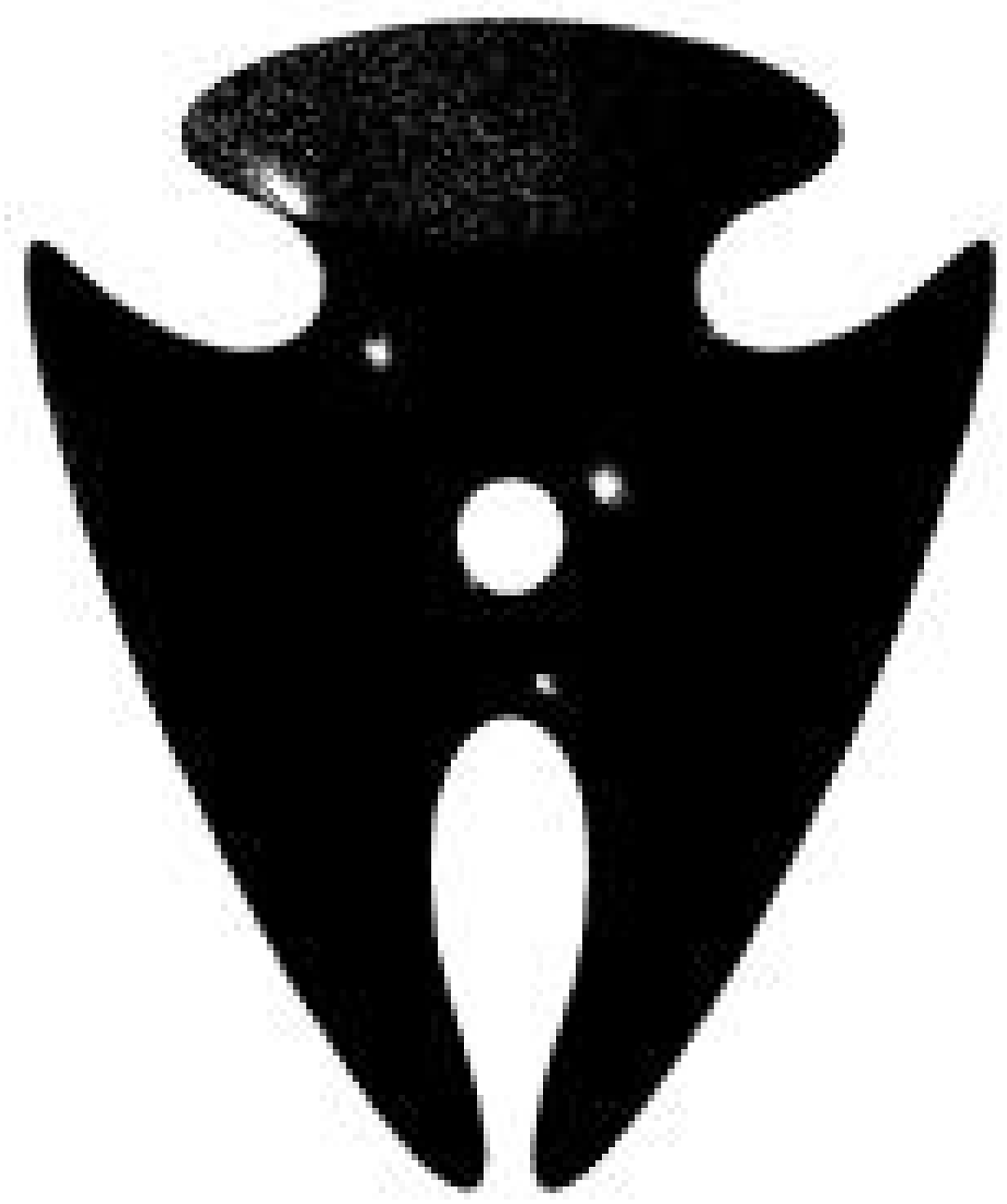}
	\hfill
	\epsfxsize = 1.5in
	\epsffile{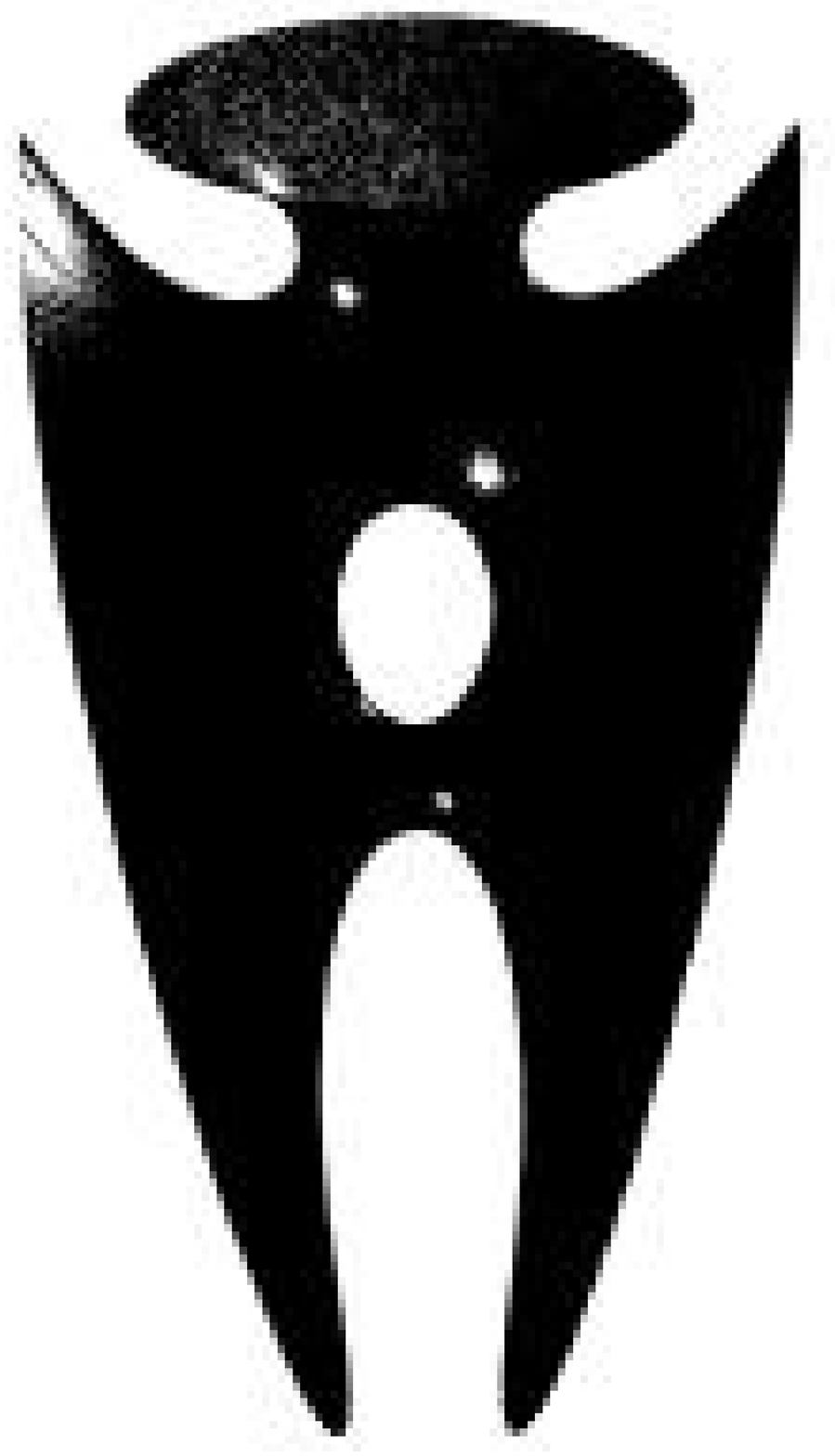}
	\hspace{0.5in}
	\caption{Less symmetric trinoids of genus 1: (left) with angles of
	     $\frac{2\pi}{3}$ between the normals of the ends, and (right) where 
	     those angles are $\frac{3\pi}{5}$, 
	     $\frac{3\pi}{5}$, and $\frac{4\pi}{5}$}  
 	\label{fig:108}
\end{figure}

To solve the period problems, we need to have the 
boundary planar 
geodesics $\alpha_1$ and $\alpha_3$ in the same plane; 
we also need to have $\alpha_2$, $\alpha_4$, and 
$\alpha_5$ all within a single plane.  This can be 
accomplished by the proper choice of $\lambda_1$, 
$\lambda_2$, and $\lambda_3$.  Using the MESH 
program \cite{MESH} and a Simplex 
algorithm, we have found values for $\lambda_i$ so 
that the Weierstrass data produces one fourth of 
the symmetric genus-1 trinoid.  Surprisingly, we 
also found one other set of values for $\lambda_i$ 
which solves the period problem.  This surface is 
not as symmetric, and suggests the existence of a 
larger family of less symmetric $n$-oids of genus 1 (see Figure 
\ref{fig:108}).  

By adjusting the value of $c$ one can still solve 
the period problem with the Simplex method and 
produce genus-1 trinoids where the normals at the 
ends to not form angles of $\frac{2\pi}{3}$ (see Figure \ref{fig:108}).

\pagebreak

\end{document}